\documentclass{elsart}

\usepackage{amssymb, amscd}
\input xy
\xyoption{all} \xyoption{poly} \xyoption{knot}

\def\p{{\mathfrak p}}

\def\CC{\mathbb{C}}   
\def\RR{\mathbb{R}}   
\def\ZZ{\mathbb{Z}}
\def\NN{\mathbb{N}}
\def\GG{\mathbb{G}}
\def\GLnC{\textbf{GL}_n \left( \mathbb{C} \right)}  
\def\GL1C{\textbf{GL}_1 \left( \mathbb{C} \right)}  
\def\GL{\textbf{GL}}
\def\BnC{\textbf{B}_n \left( \mathbb{C} \right)}    
\def\UnC{\textbf{U}_n \left( \mathbb{C} \right)}    
\def\DnC{\textbf{D}_n \left( \mathbb{C} \right)}    
\def\MnC{\textbf{M}_n \left( \mathbb{C} \right)}    
\def\p1Mx{\pi_1 \left( M, x \right)}
\def\ker{\textnormal{ker }}

\def\mc{\mathcal}

\def\gr{\mathbf{gr}}

\newcommand\im{\operatorname{im}}

\newcommand\Tor{\operatorname{Tor}}

\newcommand\Spec{\operatorname{Spec}}

\newcommand\Aut{\operatorname{Aut}}

\newtheorem{theorem}{Theorem}[section]
\newtheorem{lemma}[theorem]{Lemma}
\newtheorem{proposition}[theorem]{Proposition}
\newtheorem{corollary}[theorem]{Corollary}

\theoremstyle{definition}
\newtheorem{definition}[theorem]{Definition}

\theoremstyle{remark}

\begin{document}

\begin{frontmatter}



\title{Exponential Iterated Integrals and the Relative 
Solvable Completion of the Fundamental Group of a Manifold\thanksref{NSF}}
\thanks[NSF]{This work was partially supported by the Duke math department's
VIGRE grant, DMS-9983320 from the NSF.}
\author{Carl Miller}
\ead{carl@math.berkeley.edu}
\address{Department of Mathematics\\ 970 Evans Hall\\ University of 
California \\
Berkeley, CA 94720-3840}


\begin{abstract}
We develop a class of integrals on a manifold $M$ called {\it exponential 
iterated integrals}, an extension of K.~T.~Chen's iterated integrals.  It 
is shown that the matrix entries of any upper triangular  
representation of $\pi_1(M,x)$ can be expressed via these new integrals.  
The ring of exponential iterated integrals contains the coordinate rings 
for a class of universal representations, called the {\it relative 
solvable completions} of $\pi_1(M,x)$.  We consider exponential 
iterated integrals in the particular case of 
fibered knot complements, where the fundamental group always has a 
faithful relative solvable completion.
\end{abstract}

\end{frontmatter}

\section{Introduction}

We are concerned with using integrals to determine the
fundamental group of a smooth manifold $M$.  Let $PM$ denote the
space of piecewise differentiable paths $\lambda \colon [0,1]
\rightarrow M$.  A $1$-form $\omega \in E^1(M; \CC)$ provides a
$\CC$-valued function on $PM$ via integration:
\[
\int \omega \colon PM \rightarrow \CC
\]
\[
\lambda \mapsto \int_\lambda \omega.
\]
And $\int \omega$ induces a map on the fundamental group $\p1Mx$
if and only if $\omega$ is closed:
\[
\int \omega \colon \p1Mx \rightarrow \CC.
\]
By the de Rham theorem, closed line integrals can distinguish two
elements of $\p1Mx$ if and only if they are different in $H_1(M;
\CC)$.

K.-T. Chen improved this approach with {\it iterated
integrals} of $1$-forms (see \cite{chen2}).  For $\CC$-valued
$1$-forms $\omega_1, \omega_2, \ldots, \omega_n$,
\[
\int_\lambda \omega_1 \omega_2 \ldots \omega_n := \int_{0 \leq t_1
\leq t_2 \leq \ldots \leq t_n \leq 1} f_1 \left( t_1 \right) f_2
\left( t_2 \right) \ldots f_n \left( t_n \right) dt_1 dt_2 \ldots
dt_n
\]
where $f_i \left( t \right) dt$ is the pullback of $\omega_i$ to $E^1
\left( [0,1]; \CC \right)$ along $\lambda \colon [0,1] \rightarrow M$.  
An iterated integral is a finite sum of these expressions, with a constant
term, regarded as a function from $PM$ into $\CC$.  A {\it closed iterated
integral} is one that is constant on homotopy classes of paths $\lambda
\colon [0,1] \rightarrow M$ relative to $\{0,1\}$, and thus induces a
$\CC$-valued function on $\p1Mx$. The vector space of iterated integrals
on $PM$ is denoted by $B(M)$ and the vector space of closed iterated
integrals on the space of loops based at $x \in M$ is denoted by $H^0
\left( B (M, x) \right)$. Both are commutative Hopf algebras.

This larger class of integrals can be used to detect more
structure in $\p1Mx$ than is detected by ordinary line
integrals.  Chen proved that integration induces a Hopf algebra
isomorphism,
\begin{eqnarray}
\label{chen} H^0 \left( B (M, x) \right) \cong \mathcal{O} \left(
\mathcal{U} \left( \p1Mx \right) \right),
\end{eqnarray}
where $\mathcal{O} \left( \mathcal{U} \left( \p1Mx \right)
\right)$ denotes the coordinate ring of the unipotent completion
of $\p1Mx$.  Thus when the
representation $\p1Mx \rightarrow \mathcal{U} \left( \p1Mx
\right)$ is faithful (this occurs, for example, when $\p1Mx$ is free), 
closed iterated integrals separate the
elements of $\p1Mx$.

But there are important cases where $\p1Mx \rightarrow
\mathcal{U} \left( \p1Mx \right)$ is far from being faithful.
Indeed, if $H_1 \left(M; \ZZ \right) \cong \ZZ$ (as in the case
of knot groups), $\mathcal{U} \left( \p1Mx \right)$ is the
additive group $\mathbb{G}_a$ and the kernel of the
representation is the commutator subgroup of $\p1Mx$.  Thus closed
iterated integrals vanish on every element of the commutator
subgroup, and they provide no advantage over ordinary line
integrals in this case.

We shall overcome this limitation by considering a larger class
of integrals, called {\it exponential iterated integrals}.  The
goal is to detect a larger quotient of $\p1Mx$ than is obtained
by ordinary iterated integrals.

Exponential iterated integrals are certain infinite sums of
ordinary iterated integrals, and their properties are similar.  An
exponential iterated integral is written as
\[
\sum \int e^{\delta_1} \omega_{12} e^{\delta_2} \omega_{23}
e^{\delta_3} \ldots \omega_{(n-1)n} e^{\delta_n},
\]
where $\delta_i$ and $\omega_{i(i+1)}$ are $1$-forms.  If $L
\subset E^1 (M; \CC)$ is a $\ZZ$-module of closed $1$-forms, we
denote by $EB(M)^L$ the vector space of exponential iterated
integrals whose exponents $\delta_i$ are in $L$, and by $H^0
\left( EB(M, x)^L  \right)$ the space of exponential iterated
integrals on the loop space at $x$ that are constant on homotopy
classes of paths relative to endpoints. Each integral in $H^0
\left( EB(M, x)^L \right)$ gives a map $\pi_1(M,x) \to \CC$. We
will show that these integrals are matrix entries of solvable
representations of $\pi_1(M,x)$.

We then develop the notion of {\it relative solvable completion},
which is a special case of Deligne's relative unipotent
completion (see \cite{completions}).  Given a homomorphism $\rho
\colon G \rightarrow T$ of an abstract group $G$ into a
diagonalizable algebraic group $T$ with Zariski dense image, the
{\it solvable completion of $G$ relative to $\rho$}, denoted by
$\mathcal{S}_\rho \left( G \right)$, is the inverse limit of
algebraic representations $\phi \colon G \rightarrow S$ that fit
into a commutative diagram:
\[
\xymatrix{ & & G \ar[d]^\phi \ar[rd]^\rho \\
1 \ar[r] & U \ar[r] & S \ar[r] & T \ar[r] & 1 }
\]
where the bottom row is exact, $\phi$ has Zariski dense image,
and $U$ is a unipotent group.

The essential link between exponential iterated integrals and
relative solvable completion is given by the following theorem, a
more general version of which is proved in Section~\ref{bigthm}.

\begin{theorem}
\label{prelim} Suppose $\rho \colon \p1Mx \rightarrow T \subseteq
\left( \CC^* \right)^n$ is a diagonal algebraic representation
with Zariski dense image, and that $\delta_1, \ldots, \delta_n$
are closed $1$-forms such that
\[
\rho (\lambda) =  \left( e^{\int_\lambda \delta_1}, \ldots ,
e^{\int_\lambda \delta_n} \right) \in \left( \CC^* \right)^n.
\]
Let $L$ denote the $\ZZ$-submodule of $E^1(M; \CC)$ generated by
$\delta_1, \ldots, \delta_n$.  Then integration induces a Hopf
algebra isomorphism
\[
H^0 \left(EB(M, x)^L \right) \cong \mathcal{O} \left(
\mathcal{S}_\rho \left( \p1Mx \right) \right).
\]
\end{theorem}

In Section~\ref{radical} we consider the relationship between
unipotent completion and relative solvable completion, and prove
a result of which the following is a special case:

\begin{theorem}
If $G$ is a group such that $G/[G,G]$ is finitely generated and 
$H_1 \left( [G, G];
\CC \right)$ is finite dimensional, then there exists a
diagonalizable algebraic representation $\rho \colon G \to T$
such that $\mc{U} \left( \left[ G, G \right] \right)$ injects
into $\mathcal{S}_\rho \left( G \right)$.
\end{theorem}

All knot groups $G = \pi_1 \left( S^3 - K, x
\right)$ of tame knots $K$ satisfy the conditions of this
theorem.  The representation $\rho \colon G \to T$ can be
obtained from the Alexander module of $K$. When $K$ is a
fibered knot, $\left[ G, G \right]$ is free and thus it injects into
its unipotent completion. Hence:

\begin{corollary}
\label{here} If $K \subset S^3$ is a fibered knot, there exists a
diagonalizable algebraic representation $\rho \colon \pi_1 \left(
S^3 \smallsetminus K, x \right) \rightarrow T$ such that the 
representation
$\pi_1 \left( S^3 \smallsetminus K, x \right) \rightarrow \mathcal{S}_\rho
\left( \pi_1 \left( S^3 \smallsetminus K, x \right) \right)$ is injective.
\end{corollary}

Combining Corollary~\ref{here} with Theorem~\ref{prelim} shows that
exponential iterated integrals separate the elements of the group of a
fibered knot.  In Section~\ref{examples} we consider the example of the 
trefoil
knot (a fibered knot), providing an explicit description for the vector
space of closed exponential iterated integrals on its complement in $S^3$.

\section{Notation and Conventions}

Throughout, $M$ is a $C^\infty$-manifold with base point $x$.
$PM$ is the space of piecewise differentiable paths $\lambda
\colon [0,1] \rightarrow M$, and $P_{x,x}M$ is the loop space at
$x \in M$. $E^1(M; \CC)$ is the space of $\CC$-valued $1$-forms,
and $B^1(M; \CC)$ is the space of closed $\CC$-valued $1$-forms.
When $I$ is an integral and $\lambda \in PM$ is a path, let
$\left< I, \lambda \right>$ denote the integral of $I$ over
$\lambda$. For $y \in M$, let $\mathbf{1}_y$ denote the constant
loop at $y$. If $\delta$ is a path or $1$-form, we write
$[\delta]$ to mean the homotopy, homology, or cohomology class of
$\delta$, depending on the context.

By ``algebraic group'' we will always mean linear algebraic group
over $\CC$. We say that an algebraic group is {\it
diagonalizable} if it is isomorphic to a closed subgroup of
$\left( \CC^* \right)^n$ for some $n$. If $\textbf{G}$ is an
algebraic group then $\textbf{G}_u$ denotes the unipotent radical
of $\textbf{G}$. Let $\MnC$, $\BnC$, $\UnC$, and $\DnC$ denote,
respectively, the $n \times n$ matrix ring, upper triangular
matrix group, unipotent matrix group, and diagonal matrix group
over $\CC$.

\section{Exponential Iterated Integrals}

In this section we define exponential iterated integrals and show
how they appear as matrix entries for the transport functions of
certain trivialized vector bundles.  Using this relationship we
then prove several formal properties that will be used in later
proofs.

Throughout this section let $\lambda \colon [0,1] \to M$ be a
path.  We begin with the definition of ordinary iterated
integrals. Note that the definition given in the introduction
extends easily to iterated integrals of $1$-forms taking values
in any $\CC$-algebra.

\begin{definition}
Suppose $\omega_1, \ldots, \omega_n$ are $1$-forms on $M$ taking
values in an associative $\CC$-algebra $A$.  Let
\[
\int_\lambda \omega_1 \omega_2 \ldots \omega_n = \int_{0 \leq t_1
\leq t_2 \leq \ldots \leq t_n \leq 1} F_1 \left( t_1 \right) F_2
\left( t_2 \right) \ldots F_n \left( t_n \right) dt_1 dt_2 \ldots
dt_n,
\]
where $F_i(t) dt = \lambda^* \omega_i \in E^1\left( [0,1]; \CC
\right) \otimes A$.
\end{definition}
The expression $\int \omega_1 \ldots \omega_n$ denotes a map from $PM$ to 
$A$.  An $A$-valued iterated integral is a finite sum of these 
expressions, possibly including a constant term.
The following theorem of
Chen's provides the initial motivation for this definition
(\cite{hain}, pp. 253):
\begin{theorem}
\label{transport} Given a trivial vector bundle $\CC^n \times M
\rightarrow M$ with connection $\nabla = d -
\omega$,\footnote{This means that for a section $f \colon M
\rightarrow \CC^n$, $\nabla f = df - f \omega$.} $\omega \in E^1
\left( M; \CC \right) \otimes \MnC$, let $T \colon PM \rightarrow
GL \left( n, \CC \right)$ denote the transport function.  For any
$\lambda \in PM$, the sum
\begin{eqnarray*}
I + \int_\lambda \omega + \int_\lambda \omega \omega +
\int_\lambda \omega \omega \omega + \ldots
\end{eqnarray*}
converges absolutely, and
\begin{eqnarray}
\label{transsum} T(\lambda) = I + \int_\lambda \omega +
\int_\lambda \omega \omega + \int_\lambda \omega \omega \omega +
\ldots  \qed
\end{eqnarray}
\end{theorem}
When the matrix of $1$-forms $\omega$ is strictly upper
triangular, the series above is finite, and the transport
function is given by a matrix of ($\CC$-valued) iterated integrals. For 
example,
if
\begin{eqnarray}
\label{omeganil}
\omega = \left[
\begin{array}{cccccc}
0 & \omega_{12} & 0 & 0 & \ldots & 0 \\
0 & 0 & \omega_{23} & 0 & \ldots & 0 \\
0 & 0 & 0 & \omega_{34} & \ldots & 0 \\
\vdots & \vdots & \vdots & \vdots & \ddots & \vdots \\
0 & 0 & 0 & 0 & \ldots & \omega_{(n-1)n} \\
0 & 0 & 0 & 0 & \ldots & 0 \\
\end{array}
\right],
\end{eqnarray}
computing the series yields
\begin{eqnarray}
\label{ordseries}
T = \left[
\begin{array}{cccccc}
1 & \int \omega_{12} & \int \omega_{12} \omega_{23} & \int
\omega_{12} \omega_{23} \omega_{34} & \ldots & \int \omega_{12}
\omega_{23} \ldots \omega_{\left( n - 1 \right) n} \\
0 & 1 & \int \omega_{23} & \int \omega_{23} \omega_{34} & \ldots
& \int \omega_{23} \omega_{34} \ldots \omega_{\left( n - 1 \right) n} \\
0 & 0 & 1 & \int \omega_{34} & \ldots
& \int \omega_{34} \omega_{45} \ldots \omega_{\left( n - 1 \right) n} \\
\vdots & \vdots & \vdots & \vdots & \ddots & \vdots \\
0 & 0 & 0 & 0 & \ldots & \int \omega_{(n-1)n} \\
0 & 0 & 0 & 0 & \ldots & 1 \\
\end{array}
\right].
\end{eqnarray}
We write $B(M,x)$ for the vector space of functions $P_{x, x}M \to
\CC$ that are given by iterated integrals. The subspace of
functions that are constant on homotopy classes is denoted by
$H^0 \left( B(M, x) \right)$.\footnote{The reason for this notation is 
that $H^0 \left( B(M,x) \right)$ is the first cohomology group of the 
complex of higher iterated integrals (see \cite{chen2}, 
Section 1.5).} 

Now we can define exponential iterated integrals:

\begin{definition}
\label{expitint} For $n \geq 0$ and $\delta_1, \delta_2, \ldots,
\delta_n, \omega_{12}, \ldots, \omega_{\left( n - 1
\right) n} \in E^1 \left( M; \CC \right)$,
\begin{eqnarray}
\nonumber \lefteqn{  \int_\lambda e^{\delta_1} \omega_{12}
e^{\delta_2} \omega_{23} e^{\delta_3} \ldots e^{\delta_{n-1}}
\omega_{\left( n
- 1 \right) n} e^{\delta_n}} \\
\label{singexp} & := & \sum_{m_1, \ldots, m_n \geq 0}
\int_\lambda \underbrace{\delta_1 \delta_1 \ldots \delta_1}_{m_1
\textnormal{ terms}} \omega_{12} \underbrace{\delta_2 \delta_2
\ldots \delta_2}_{m_2 \textnormal{ terms}} \omega_{23} \ldots
\omega_{\left( n - 1 \right) n} \underbrace{\delta_n \delta_n
\ldots \delta_n}_{m_n \textnormal{ terms}}.
\end{eqnarray}
\end{definition}
An exponential iterated integral is a finite sum of these
expressions, regarded as a function from $PM$ to $\CC$.  By the
{\it length} of an exponential iterated integral we mean the
number of linear $1$-forms in its longest term.  The integral
above has length $n-1$.\footnote{There is a possible ambiguity
because two different integral expressions can compute the same
map $PM \rightarrow \CC$.  So we will say that the length of an
exponential iterated integral is the minimum length of its
literal expressions.}

To see where Definition~\ref{expitint} comes from, suppose $E =
\CC^n \times M \rightarrow M$ is a trivial bundle with connection
$\nabla = d - \omega$ where $\omega$ is a superdiagonal matrix:
\begin{eqnarray}
\label{omegaform}
\omega = \left[ \begin{array}{cccccc}
\delta_1 & \omega_{12} & 0 & 0 & \cdots & 0 \\
0 & \delta_2 & \omega_{23} & 0 & \cdots & 0 \\
0 & 0 & \delta_3 & \omega_{34} & \cdots & 0 \\
0 & 0 & 0 & \delta_4 & \cdots & 0 \\
\vdots & \vdots & \vdots & \vdots & \ddots & \vdots \\
0 & 0 & 0 & 0 & \cdots & \delta_n \\
\end{array} \right].
\end{eqnarray}
Computing the matrix entries of the series (\ref{transsum}) yields
\begin{eqnarray}
\label{transform}
T = \left[
\begin{array}{cccccc}
\int e^{\delta_1} & \int e^{\delta_1} \omega_{12} e^{\delta_2} &
\int e^{\delta_1} \omega_{12} e^{\delta_2} \omega_{23}
e^{\delta_3} & \ldots & \int e^{\delta_1} \omega_{12}
\ldots \omega_{\left(n-1 \right) n} e^{\delta_n} \\
0 & \int e^{\delta_2} & \int e^{\delta_2} \omega_{23} e^{\delta_3}
& \ldots & \int e^{\delta_2} \omega_{23} \ldots
\omega_{\left(n-1 \right) n} e^{\delta_n} \\
0 & 0 & \int e^{\delta_3} & \ldots & \int e^{\delta_3} \omega_{34}
\ldots \omega_{\left(n-1 \right) n} e^{\delta_n} \\
\vdots & \vdots & \vdots & \ddots & \vdots \\
0 & 0 & 0 & \cdots & \int e^{\delta_n} \\
\end{array} \right].
\end{eqnarray}
Right away, then, we know from Theorem~\ref{transport} that exponential
iterated integrals are well-defined.
\begin{proposition}
For any 1-forms $\{ \delta_j \}_j, \{ \omega_{j(j+1)} \}_j
\subseteq E^1 \left( M; \CC \right)$, the sum (\ref{singexp}) converges
absolutely. \qed
\end{proposition}
Also, since transport is invariant under reparametrization of paths, we
have
\begin{proposition}
\label{indep} For any $1$-forms $\delta_1, \ldots, \delta_n,
\omega_{12}, \ldots, \omega_{(n-1)n}$, the integral
\[
\int_\lambda e^{\delta_1} \omega_{12} \ldots \omega_{(n-1)n}
e^{\delta_n}
\]
is independent of the parametrization of $\lambda$. \qed
\end{proposition}

And if $\lambda(0) = y$, then $T\left( \lambda \lambda^{-1}
\right) = T \left( \mathbf{1}_y \right)$; thus,
\begin{proposition}
\label{inverses} For any $1$-forms $\delta_1, \ldots, \delta_n,
\omega_{12}, \ldots, \omega_{(n-1)n}$,
\[
\int_{\lambda \lambda^{-1}} e^{\delta_1} \omega_{12} \ldots
\omega_{(n-1)n} e^{\delta_n} = \int_{\mathbf{1}_{\lambda(0)}}
e^{\delta_1} \omega_{12} \ldots \omega_{(n-1)n} e^{\delta_n}. \qed
\]
\end{proposition}

A formula for the integral $\int e^{\delta_1} \omega_{12} \ldots
\omega_{(n-1)n} e^{\delta_n}$ over $\lambda^{-1}$ is easily
verified from the definition:

\begin{proposition}
\label{antipode} For any $1$-forms $\delta_1, \ldots, \delta_n,
\omega_{12}, \ldots, \omega_{(n-1)n} \in E^1 \left( M; \CC \right)$,
\[
\int_{\lambda^{-1}} e^{\delta_1} \omega_{12} \ldots
\omega_{(n-1)n} e^{\delta_n} = \int_\lambda e^{-\delta_n} \left(
-\omega_{(n-1)n} \right) \ldots \left( - \omega_{12} \right)
e^{-\delta_1}. \qed
\]
\end{proposition}

For the remaining propositions let
\[
I = \int e^{\delta_1} \omega_{12} \ldots \omega_{(n-1)n}
e^{\delta_n}.
\]

The expression (\ref{transform}) for parallel transport in $E$ can
be applied to prove a formula for the integral of $I$ over a
concatenation of paths:

\begin{proposition}
\label{comult} For any paths $\alpha, \beta \in PM$ with
$\alpha(1) = \beta(0)$,
\[
\left< I, \alpha \beta \right> = \sum_{k = 1}^n \int_{\alpha}
e^{\delta_1} \omega_{12} \ldots \omega_{\left( k - 1 \right) k}
e^{\delta_k} \int_{\beta} e^{\delta_k} \omega_{k\left(k+1\right)}
\ldots \omega_{\left( n - 1 \right) n} e^{\delta_n}.
\]
\end{proposition}

\begin{pf}
The transport map satisfies $T \left( \alpha
\beta \right) = T \left(
\alpha \right) T \left( \beta \right)$. The above formula comes from comparing the
upper-right-hand matrix entries of $T \left( \alpha \beta \right)$ and $T
\left( \alpha \right) T \left( \beta \right)$.
\end{pf}

Given $s,t \in \left[ 0, 1 \right]$, let $\lambda_s^t$ denote the
subpath of $\lambda$ from $\lambda(s)$ to $\lambda(t)$, defined by
\[
\lambda_s^t(u) = \lambda(s + (t-s)u).
\]
The concatenation of $\lambda_0^t$ and $\lambda_t^1$ is equal to
$\lambda$ after reparametrization, thus

\begin{corollary}
\label{cat} For any $t_0 \in [0,1]$,
\[
\left< I, \lambda \right> = \sum_{i=1}^n \int_{\lambda_0^{t_0}}
e^{\delta_1} \omega_{12} \ldots \omega_{(i-1)i} e^{\delta_i}
\int_{\lambda_{t_0}^1} e^{\delta_i} \omega_{i(i+1)} \ldots
\omega_{(n-1)n} e^{\delta_n}. \qed
\]
\end{corollary}

A similar proposition (useful in doing induction on these integrals by
length) follows immediately from the definition:

\begin{proposition}
\label{split} Let $f_{(i-1)i}(t)dt = \lambda^* \omega_{(i-1)i}$.
\[
\left< I, \lambda \right> = \int_0^1 \left( \int_{\lambda_0^t}
e^{\delta_1} \omega_{12} \ldots e^{\delta_{i-1}} \right)
f_{(i-1)i}(t) \left( \int_{\lambda_t^1} e^{\delta_i} \ldots
\omega_{(s-1)s} e^{\delta_s} \right) dt. \qed
\]
\end{proposition}

Note that
\begin{eqnarray*}
\int_\lambda e^\delta & = & \sum_{n \in \NN} \int_{0 \leq t_1 \leq
t_2 \leq \ldots \leq t_n \leq 1} f \left( t_1 \right) f \left( t_2
\right) \ldots f \left( t_n \right) dt_1 dt_2 \ldots dt_n \\
& = & \sum_{n \in \NN} \frac{1}{n!} \int_{t_1, t_2, \ldots, t_n
\in \left[ 0, 1 \right]} f \left( t_1 \right) f \left( t_2 \right)
\ldots f
\left( t_n \right) dt_1 dt_2 \ldots dt_n \\
& = & \sum_{n \in \NN} \frac{1}{n!} \left( \int_{0 \leq t \leq 1}
f \left( t \right) dt \right)^n \\
& = & e^{\int_\lambda \delta} \\
\end{eqnarray*}
(where $f(t)dt = \lambda^* \delta$), hence our notation.

\begin{proposition}
\label{zerolength} For any $\delta \in E^1 (M; \CC)$, $\lambda \in PM$,
\[
\int_{\lambda} e^\delta = e^{\int_\lambda \delta}. \qed
\]
\end{proposition}

Using this fact we can simplify integrals that have an exact
exponent:

\begin{proposition}
\label{expexact} Suppose $g \colon M \rightarrow \CC$ is a
$C^\infty$ function.  Then
\begin{eqnarray*}
\lefteqn{ \int e^{\delta_1} \omega_{12} \ldots \omega_{(i-1)i}
e^{dg} \omega_{i(i+1)} \ldots \omega_{(n-1)n} e^{\delta_n}} & & \\
& = & \int e^{\delta_1} \omega_{12} \ldots \left( e^{-g}
\omega_{(i-1)i} \right) \left( e^g \omega_{i(i+1)} \right) \ldots
\omega_{(n-1)n} e^{\delta_n}.
\end{eqnarray*}
\end{proposition}

\begin{pf}
Again let $f_{(i-1)i}(t)dt = \lambda^* \omega_{(i-1)i}$,
$f_{i(i+1)}(t)dt = \lambda^* \omega_{i(i+1)}$.  Applying
Proposition~\ref{split} twice,
\begin{eqnarray*}
\lefteqn{ \int_\lambda e^{\delta_1} \omega_{12} \ldots
\omega_{(i-1)i} e^{dg} \omega_{i(i+1)} \ldots \omega_{(n-1)n}
e^{\delta_n} = }& &  \\
& & \int_{0 \leq t \leq t' \leq 1} \left( \int_{\lambda_0^t}
e^{\delta_1} \omega_{12} \ldots e^{\delta_{i-1}} \right)
f_{(i-1)i}(t) d t  \\
& & \cdot \left( \int_{\lambda_t^{t'}} e^{df} \right)
f_{i(i+1)}(t') d t' \left( \int_{\lambda_{t'}^1} e^{\delta_{i+1}}
\ldots \omega_{(n-1)n} e^{\delta_n} \right)\\
 & = & \int_{0 \leq
t \leq t' \leq 1} \left( \int_{\lambda_0^t} e^{\delta_1}
\omega_{12} \ldots e^{\delta_{i-1}} \right)
f_{(i-1)i}(t) d t \\
& & \cdot \textnormal{ } e^{- g\left(t\right) + g\left(t'\right)}
f_{i(i+1)}(t') d t' \left( \int_{\lambda_{t'}^1} e^{\delta_{i+1}}
\ldots \omega_{(n-1)n} e^{\delta_n}
\right) \\
& = &\int_\lambda e^{\delta_1} \omega_{12} \ldots \left( e^{-g}
\omega_{(i-1)i} \right) \left( e^g \omega_{i(i+1)} \right) \ldots
\omega_{(n-1)n} e^{\delta_n}, \\
\end{eqnarray*}
as desired.
\end{pf}

Lastly, we show how the computation (\ref{transform}) of the
transport function in terms of exponential iterated integrals can
be generalized. Suppose $\nabla = d - \omega$ is a connection on
the trivial bundle $\CC^n \times M \to M$, where $\omega$ is an
upper triangular matrix of $1$-forms:
\[
\omega = \left[ \begin{array}{cccccc} \omega_{11} & \omega_{12} &
\omega_{13} & \cdots & \omega_{1n} \\
0 & \omega_{22} & \omega_{23} & \cdots & \omega_{2n} \\
0 & 0 & \omega_{33} & \cdots & \omega_{3n} \\
\vdots & \vdots & \vdots & \ddots & \vdots \\
0 & 0 & 0 & \cdots & \omega_{nn} \\
\end{array} \right].
\]
Matrix multiplication shows that the $(i,j)$th entry of
\[
\int \omega^m := \int \underbrace{\omega \omega \ldots \omega}_{m
\textrm{ terms}}
\]
is
\[
\sum_{i = k_1 \leq \ldots \leq k_{m+1} = j} \int \omega_{k_1k_2}
\omega_{k_2k_3} \ldots \omega_{k_{m}k_{m+1}}.
\]
(We take the sum to be zero if $i > j$.)  Thus by
Theorem~\ref{transport},
\[
T = I + \left( \sum_{m > 0,  i = k_1 \leq \ldots \leq k_{m+1} =
j} \int \omega_{k_1k_2} \omega_{k_2k_3} \ldots \omega_{k_mk_{m+1}}
\right)_{1 \leq i, j \leq n}.
\]
Grouping repeated terms,
\begin{eqnarray*}
T & = & \left( \sum_{p > 0, i = k_1 < \ldots < k_p = j, r_1, \ldots, r_p 
\geq 0} \int \omega_{k_1k_1}^{r_1} \omega_{k_1k_2}
\omega_{k_2k_2}^{r_2} \omega_{k_2k_3} \ldots \omega_{k_{p-1}k_p}
\omega_{k_pk_p}^{r_p}
\right)_{1 \leq i, j \leq n} \\
& = & \left( \sum_{p > 0, i = k_1 < \ldots < k_p = j} \int
e^{\omega_{k_1k_1}} \omega_{k_1k_2} e^{\omega_{k_2k_2}}
\omega_{k_2k_3} \ldots \omega_{k_{p-1}k_p} e^{\omega_{k_pk_p}}
\right)_{1 \leq i, j \leq n}. \\
\end{eqnarray*}
The sums in this last expression are finite, so we have expressed
$T$ as a matrix of exponential iterated integrals.

We state this result for future reference:
\begin{proposition}
\label{upper} Suppose $\nabla = d - \omega$ is a connection on
the trivial bundle $\CC^n \times M \to M$, where $\omega =
(\omega_{ij})_{1 \leq i,j \leq n}$ is an upper triangular matrix
of $1$-forms.  Then the transport function $T \colon PM \to \GLnC$
is equal to an upper triangular matrix of exponential iterated
integrals whose exponents are from the set $\left\{\omega_{11},
\omega_{22}, \ldots, \omega_{nn} \right\}$. $\qed$
\end{proposition}

\section{Relative Solvable Representations}

\label{relsolv}

In order to describe the set of maps $\pi_1(M,x) \to \CC$ that are given
by exponential iterated integrals, we need first to explore the properties
of a particular class of algebraic representations associated to a
discrete group.

Suppose $G$ is a group and $\rho \colon G \to (\CC^*
)^n$ is a diagonal representation. Let $T \subseteq (\CC^* )^n$
denote the Zariski closure of $\rho(G)$.  By a {\it solvable
representation relative to} $\rho$ we mean an algebraic
representation $G \to S$ that fits into a commutative diagram
\begin{eqnarray}
\label{relsolrep}
\xymatrix{ & & G \ar[rd] \ar[d] \\
 1 \ar[r] & U \ar[r] & S \ar[r] & T \ar[r] & 1 }
\end{eqnarray}
where the bottom row is exact, $U$ is a unipotent group, and the
image of $G$ in $S$ is Zariski dense.

The canonical example of a relative solvable representation occurs
when $G$ has an upper triangular action on $\CC^n$. The quotient
of $\BnC$ by $\UnC$ is isomorphic to $\DnC$, hence any
representation $G \to \BnC$ fits into a diagram
\[
\xymatrix{ & & G \ar[rd] \ar[d] \\
 1 \ar[r] & \UnC \ar[r] & \BnC \ar[r] & \DnC \ar[r] & 1.}
\]
The Zariski closure of $G$ in $\BnC$ is a relative solvable
representation.  In fact every relative solvable representation
can be obtained in this manner, as the next proposition and
corollary show.

\begin{proposition}
Suppose $1 \to U \to S \to T \to 1$ is an exact sequence of
algebraic groups with $U$ unipotent and $T$ diagonalizable.  Then
any linear representation $S \to \GL(V)$ has a common eigenvector.
\end{proposition}

\begin{pf}
Suppose $S \to \GL(V)$ is a representation.  Since $U$ is
unipotent we may find a nonzero vector $v \in V$ fixed by $U$.
Let $W \subseteq V$ be the subspace spanned by $\{ xv \mid x \in S
\}$. Note that $U$ fixes $W$, since for any $x \in S$, $y \in U$,
\[
yxv = x(x^{-1}yx)v = xv.
\]
So the action of $S$ on $W$ factors through $T$; and since $T$ is
diagonalizable it has a common eigenvector in $W$.
\end{pf}

\begin{corollary}
\label{relsolisupper} Suppose
\begin{eqnarray}
\xymatrix{ & & G \ar[rd] \ar[d] \\
 1 \ar[r] & U \ar[r] & S \ar[r] & T \ar[r] & 1 }
\end{eqnarray}
is a relative solvable representation.  Then there exists an
embedding $S \hookrightarrow \BnC$ for some $n$, under which $S
\cap \UnC = U$.
\end{corollary}

\begin{pf}
Let $S \hookrightarrow \GL(V)$ be a faithful representation.
Proceeding by induction with the proposition we may find a basis
$\{ v_1, \ldots, v_n \} \subseteq V$ that makes the action of $S$
upper triangular.  Thus we obtain the embedding $S \hookrightarrow
\BnC$. The intersection of $S$ with $\UnC$ is the unipotent
radical of $S$, which is precisely $U$.
\end{pf}

We note in passing that the discussion above also explains the
term ``relative solvable.''  Since every connected solvable group
is isomorphic to a closed subgroup of $\BnC$ (see \cite{humph},
Section 19), every connected solvable representation with Zariski
dense image may be expressed as a relative solvable
representation. (The converse is not true, since a relative
solvable representation need not be connected.)

By analogy we use the term {\it relative prosolvable
representation} to mean a commutative diagram of the form of
(\ref{relsolrep}) with $U$ prounipotent and $T$ diagonalizable.

\begin{proposition}
\label{itsplits} Suppose $1 \to \mc{U} \to \mc{S} \to T \to 1$ is
an exact sequence with $\mc{U}$ prounipotent and $T$
diagonalizable.  Then the sequence splits, and any two splittings
are conjugate by an element of $\mc{U}$.
\end{proposition}

\begin{pf}
The assertion is well-known when $\mc{U}$ is unipotent (see
\cite{borelserre}, Prop. 5.1).  Choose compatible splittings for
the unipotent quotients of $\mc{U}$ and take the inverse limit.
\end{pf}

Suppose $\rho \colon G \to T \subseteq \left( \CC^* \right)^n$ is
a diagonalizable representation as before. There is a unique
relative prosolvable representation $G \to \mathcal{S}_\rho(G) \to
T$ satisfying the following universal mapping property: if
$\mc{S}$ is any prosolvable representation relative to $\rho$,
then there is a homomorphism $\mc{S}_\rho(G) \to \mc{S}$ of
$G$-representations that makes the following diagram commute:
\[
\xymatrix{ \mc{S}_\rho(G) \ar[d] \ar[dr] & \\
\mc{S} \ar[r] & T.}
\]
$\mathcal{S}_\rho(G)$ is the inverse limit of all solvable
representations relative to $\rho$. This is a special case of
{\it relative Malcev completion}, a notion of Deligne's; see
\cite{completions}, Section 2 for a full development.

As an initial example of relative prosolvable completion, suppose
$G$ is finitely generated and abelian.  Then $\mathcal{S}_\rho(G)$
must be abelian, so the splitting of Proposition~\ref{itsplits}
is a direct product: $\mathcal{S}_\rho(G) \cong \mc{U} \times T$
where $\mc{U}$ is prounipotent.  $\mc{U}$ is abelian and therefore
additive, and its dimension cannot exceed the free rank of $G$.
This proves the next proposition.
\begin{proposition}
If $\rho \colon G \to T$ is a homomorphism from a finitely
generated abelian group into a diagonalizable group with Zariski
dense image, then $\mathcal{S}_\rho(G) \cong \mathbb{G}_a^m
\times T$ where $m$ is the free rank of $G$. $\qed$
\end{proposition}

\section{Algebraic Properties of Exponential Iterated Integrals}

\label{algprop}

Suppose $L \subset E^1 (M; \CC)$ is a $\ZZ$-module of $1$-forms. Let
$EB(M)^L$ denote the vector space of exponential iterated integrals with
exponents from $L$. A closed exponential iterated integral is one that is
constant on homotopy classes of paths $\lambda \colon \left[ 0, 1 \right]
\rightarrow M$ relative to $\{ 0, 1 \}$.  Let $H^0 \left( EB(M)^L \right)$
denote the subspace of closed iterated integrals in
$EB(M)^L$.
In this section we show that
$EB(M)^L$ and $H^0 \left( EB(M)^L \right)$ are Hopf algebras, and that the
groups $\Spec EB(M)^L$ and $\Spec H^0 \left( EB(M)^L \right)$ are each an
extension of a prodiagonalizable group by a prounipotent group.

We also let $EB(M,x)^L$ denote the space of functions on the loop space
$P_{x,x}M$ given by exponential iterated integrals with exponents from
$L$.  ($EB(M,x)^L$ is the quotient of $EB(M)^L$ by integrals such as $\int
e^{df}$ that vanish on every loop at $x$.  We will also refer to the
elements of $EB(M,x)^L$ as exponential iterated integrals.)  Write
$H^0\left( EB(M,x)^L \right)$ for the subspace of functions constant on
homotopy classes in $P_{x,x}M$.  For convenience we will omit basepoints
in the rest of this section, but everything we say for $EB(M)^L$ and
$H^0\left( EB(M)^L\right)$ applies as well to $EB(M,x)^L$ and $H^0\left(
EB(M,x)^L \right)$.

\begin{lemma}
For any $\ZZ$-module $L \subseteq E^1 \left( M ; \CC \right)$,
$EB\left(M \right)^L$ is closed under (pointwise) multiplication.
\end{lemma}

\begin{pf}
It is sufficient to show that for any $\delta_i, \delta'_i \in
L$, $\omega_{j(j+1)}, \omega'_{j(j+1)} \in E^1 \left( M; \CC \right)$,
\[
\int e^{\delta_1} \omega_{12} \ldots \omega_{\left(n-1\right)n}
e^{\delta_n} \int e^{\delta'_1} \omega'_{12} \ldots
\omega'_{\left(n'-1\right)s'} e^{\delta'_{n'}} \in EB\left(M
\right)^L
\]
A formula for this product can be written down, but it is rather
cumbersome and unnecessary for our purposes.  We prove the result
by induction on $n+n'$.  Let $\lambda \colon [0,1] \to M$ be a
path.

Note first that
\[
\int_\lambda e^\delta \int_\lambda e^{\delta'} = e^{\int_\lambda
\delta} e^{\int_\lambda \delta'} = e^{\int_\lambda \delta +
\delta'} = \int_\lambda e^{\delta + \delta'}.
\]
This proves the base case.  For $n+n' \geq 1$ we apply
Corollary~\ref{cat} and Proposition~\ref{split} to split the
product into a sum of products of smaller-length integrals. Assume
without loss of generality that $n \geq 1$.  Let $f_{(n-1)n}(t)
dt = \lambda^* \omega_{(n-1)n}$.
\begin{eqnarray*}
\lefteqn{\int_\lambda e^{\delta_1} \omega_{12} \ldots
\omega_{\left(n-1\right)n} e^{\delta_n} \int_\lambda e^{\delta'_1}
\omega'_{12} \ldots \omega'_{\left(n'-1\right)n'} e^{\delta'_{n'}}} & & \\
& = & \int_{0 \leq t \leq 1}  \left( \int_{\lambda_0^t}
e^{\delta_1} \omega_{12} \ldots e^{\delta_{n-1}} \right)
f_{(n-1)n}(t) dt \left(
\int_{\lambda_t^1} e^{\delta_n} \right) \\
& & \cdot \left( \sum_{i=1}^{n'} \int_{\lambda_0^t} e^{\delta'_1}
\omega'_{12} \ldots e^{\delta'_i} \int_{\lambda_t^1}
e^{\delta'_i} \omega'_{i(i+1)}
\ldots e^{\delta'_{n'}} \right)\\
& = & \sum_{i=1}^{n'} \int_{0 \leq t \leq 1} \left(
\int_{\lambda_0^t} e^{\delta_1} \omega_{12} \ldots
e^{\delta_{n-1}} \right) \left( \int_{\lambda_0^t}
e^{\delta'_1} \omega'_{12} \ldots e^{\delta'_i} \right) \\
& & \cdot f_{(n-1)n}(t)dt \left( \int_{\lambda_t^1} e^{\delta_n}
\right) \left( \int_{\lambda_t^1} e^{\delta'_i} \omega'_{i(i+1)}
\ldots e^{\delta'_{n'}} \right) \\
\end{eqnarray*}
By inductive assumption both of
\[
\int e^{\delta_1} \omega_{12} \ldots e^{\delta_{n-1}} \int
e^{\delta'_1} \omega'_{12} \ldots e^{\delta'_i} \textnormal{ and
} \int e^{\delta_n} \int e^{\delta'_i} \omega'_{i(i+1)} \ldots
e^{\delta'_{n'}}
\]
may be expressed as exponential iterated integrals from
$EB(M)^L$. Applying Proposition~\ref{split} then transforms each
summand in the resulting expression into an exponential iterated
integral.
\end{pf}

Therefore $EB (M)^L$ is a $\CC$-algebra.  (A homomorphism $\CC \to
EB(M)^L$ is given by $z \mapsto z \int e^0$.)  Clearly the
product of two closed exponential iterated integrals is closed,
thus $H^0 \left( EB (M)^L \right)$ is likewise a $\CC$-algebra.

The definition for comultiplication $\Delta \colon EB(M)^L
\rightarrow EB(M)^L \otimes EB(M)^L$ comes from
Proposition~\ref{comult}:
\begin{eqnarray}
\nonumber \lefteqn{ \Delta \int e^{\delta_1} \omega_{12} \ldots
\omega_{\left( n - 1 \right) n} e^{\delta_n} := }& & \\
\label{comultdef}
& & \sum_{i = 1}^n \int e^{\delta_1} \omega_{12} \ldots
\omega_{\left( i - 1 \right) i} e^{\delta_i} \otimes \int
e^{\delta_i} \omega_{i\left(i+1\right)} \ldots \omega_{\left( n -
1 \right) n} e^{\delta_n}.
\end{eqnarray}
For $I \in EB(M)^L$ and $\alpha, \beta, \gamma \in PM$, we have
\[
\left< \Delta I, (\alpha, \beta) \right> = \left< I, \alpha \beta
\right>.
\]
The coassociative property of $\Delta$ follows:
\[
\left< \left( \Delta \otimes 1 \right) \Delta I, (\alpha, \beta,
\gamma) \right> = \left< I, \alpha \beta \gamma \right> = \left<
\left( 1 \otimes \Delta \right) \Delta I, (\alpha, \beta, \gamma)
\right>.
\]
And if $I$ is constant on homotopy classes in $PM$, the function
\[
(\alpha, \beta) \mapsto \left< \Delta I, (\alpha, \beta ) \right>
= \left< I, \alpha \beta \right>
\]
is constant on homotopy classes in $PM \times PM$; thus $\Delta$
restricts to a comultiplication
\[
\Delta \colon H^0 \left( EB \left( M \right)^{L} \right)
\rightarrow H^0 \left( EB \left( M \right)^{L}  \right) \otimes
H^0 \left( EB \left( M \right)^{L}  \right).
\]
The constant map on $EB(M)^L$ and $H^0\left( EB(M)^L \right)$ is
given by $I \mapsto \left< I, \mathbf{1}_x \right>$.  (Note that
$\left< I, \mathbf{1}_y \right> = \left< I, \mathbf{1}_x \right>$
for any $y \in M$.) Proposition~\ref{antipode} gives an antipode
map on $EB \left( M \right)^{L}$:
\[
i \colon \int e^{\delta_1} \omega_{12} \ldots \omega_{(n-1)n}
e^{\delta_n} = (-1)^{n-1} \int e^{-\delta_n} \omega_{(n-1)n}
\ldots \omega_{12} e^{-\delta_1}.
\]
By Proposition~\ref{inverses},
\[
\left< \left( 1 \otimes i \right) \Delta I , \lambda \right> =
\left< I, \lambda \lambda^{-1} \right> = \left< I, \mathbf{1}_x
\right>.
\]
And as with $\Delta$, $i$ restricts to an antipode map on $H^0 \left( 
EB(M)^L \right)$. We have thus proven:
\begin{proposition}
For any $\ZZ$-module of $1$-forms $L \subseteq E^1(M; \CC)$,
\[
\left(
EB(M)^L, \Delta, i \right) \textnormal{ and } \left( H^0 \left( EB(M)^L
\right), \Delta, i \right)
\]
are Hopf algebras. \qed
\end{proposition}

To demonstrate the structure of the group $\Spec EB(M)^L$ we need
a few definitions. Let $\mathcal{E}(M)^L \subseteq EB(M)^L$ denote
the $\CC$-algebra generated by integrals of the form $\int
e^\delta$ with $\delta \in L$.  Since $\Delta \int e^\delta =
\int e^\delta \otimes \int e^\delta$, $\mathcal{E}(M)^L$ is the
coordinate ring of a prodiagonalizable group.  The inclusion $i
\colon \mathcal{E}(M)^L \to EB(M)^L$ gives a surjective
homomorphism
\[
i^* \colon \Spec EB(M)^L \to \Spec \mathcal{E}(M)^L.
\]
We show that the kernel of this homomorphism is a prounipotent
group.  The coordinate ring of the kernel is $EB(M)^L/I_0$, where
$I_0$ is the ideal generated by integrals of the form $\int
e^\delta - 1$, $\delta \in L$. Consider the filtration of ideals
\[
I_0 \subseteq I_1 \subseteq I_2 \subseteq \ldots,
\]
where
\[
I_j = \left< I - \left< I, \mathbf{1}_x \right> \mid I \in EB(M)^L
\textnormal{ of length} \leq j \right>
\]
This filtration corresponds to a filtration of subgroups in
$\Spec EB(M)^L/I_0$.  For $j \geq 1$, the action of $\Delta$ on
$I_j / I_{j-1}$ is given by
\begin{eqnarray*}
\lefteqn{\Delta (\int e^{\delta_1} \omega_{12} \ldots
e^{\delta_{j+1}}
+ I_{j-1})} \\
& = & \int e^{\delta_1} \omega_{12} \ldots e^{\delta_{j+1}}
\otimes 1 + 1 \otimes \int e^{\delta_1} \omega_{12} \ldots
e^{\delta_{j+1}} + (I_{j-1} \otimes I_j + I_j \otimes I_{j-1}).
\end{eqnarray*}
Thus $EB(M)^L/I_0$ is indeed a prounipotent group.  We have
constructed an exact sequence of groups
\[
1 \to \Spec EB(M)^L/I_0 \to \Spec EB(M)^L \to \Spec
\mathcal{E}(M)^L \to 1
\]
where $\Spec EB(M)^L/I_0$ is prounipotent and $\Spec
\mathcal{E}(M)^L$ is prodiagonalizable.

Now if we assume that $L$ consists of closed $1$-forms, the
algebra $\mathcal{E}(M)^L$ is contained in $H^0 \left( EB(M)^L
\right)$, and we obtain a similar exact sequence for $\Spec
H^0\left( EB(M)^L\right)$. Let $I'_0 = I_0 \cap H^0 \left(
EB(M)^L \right)$.  The sequence
\[
1 \to \Spec H^0\left( EB(M)^L \right)/I'_0 \to \Spec H^0 \left(
EB(M)^L \right) \to \Spec \mathcal{E}(M)^L \to 1
\]
is a quotient of the one above by a prounipotent subgroup.

We summarize this discussion.
\begin{theorem}
\label{prosolvable} For any $\ZZ$-module $L \subseteq E^1(M; \CC)$, the
group $\Spec \mathcal{E}(M)^L$ is prodiagonalizable, and there is
an epimorphism
\[
\Spec EB(M)^L \to \Spec \mathcal{E}(M)^L
\]
whose kernel is prounipotent. If $L$ consists of closed
$1$-forms, there is an epimorphism
\[
\Spec H^0\left( EB(M)^L \right) \to \Spec \mathcal{E}(M)^L
\]
whose kernel is prounipotent. $\qed$
\end{theorem}

\section{The Solvable de Rham Theorem}

\label{bigthm}

Suppose that $L$ is a $\ZZ$-module of closed $1$-forms. Continuing
the notation of Section~\ref{algprop}, let $\mathcal{E}(M,x)^L
\subseteq H^0 \left( EB(M, x)^L \right)$ denote the subalgebra
generated by the integrals $\int e^\delta, \delta \in L$.
Consider the homomorphism
\[
\rho \colon \p1Mx \to \Spec \mathcal{E}(M,x)^L
\]
where $[\lambda]$ maps to the ideal of integrals that vanish on
$\lambda$. This representation is prodiagonalizable and has
Zariski dense image.  If there is a finite set $\{ \delta_1, \ldots,
\delta_n \} \subseteq L$ that spans the image of $L$ in $H^1( M ; \CC)$
then $\rho$ is an algebraic representation, and it may be expressed as
$\rho \colon \p1Mx \to T \subseteq \left( \CC^* \right)^n$,
\[
\rho(\lambda) = \left( \int_\lambda e^{\delta_1}, \ldots,
\int_\lambda e^{\delta_n} \right).
\]
We shall call $\rho$ ``the representation defined by $L$.''

As noted in the introduction, Chen proved that the Hopf algebra
of closed iterated integrals on $P_{x,x}M$ is the coordinate ring
of the unipotent completion of $\p1Mx$.  Previous work has been
done on extending this isomorphism: R.~Hain in \cite{malcev}
constructed a class of integrals that compute the coordinate ring
of the Malcev completion of $\p1Mx$ relative to any algebraic
representation $\rho \colon \p1Mx \to S$.  These integrals are written in
the form
\[
\int \left( \omega_1 \omega_2 \ldots \omega_r \mid \phi \right),
\]
where the $\omega_i$ are $1$-forms on a principal $S$-bundle over
$M$, and $\phi$ is a matrix entry of $S$.

The solvable de Rham theorem shows that the coordinate ring of the
Malcev completion may be computed using the more geometric (and
more manageable) class of exponential iterated integrals, in the
case where $\rho$ is a diagonalizable representation defined by
some $\ZZ$-module $L \subseteq E^1(M; \CC)$.

\begin{theorem}[The $\pi_1$ solvable de Rham theorem]
\label{mainthm} Suppose $L$ is a $\ZZ$-module of closed $1$-forms
whose image in $H^1(M; \CC)$ is finitely generated, and that $\rho
\colon \p1Mx \rightarrow T \subseteq \left( \CC^* \right)^n$ is
the representation defined by $L$.  Then integration induces a
Hopf algebra isomorphism
\[ H^0 \left( EB
\left( M,x \right)^{L} \right) \cong  \mathcal{O} \left(
\mathcal{S}_\rho  \left( \pi_1(M,x) \right) \right).
\]
\end{theorem}

\begin{pf}
The representation
\[
\pi_1(M,x) \to \Spec H^0 \left( EB(M,x)^L \right)
\]
has Zariski dense image, and by Theorem~\ref{prosolvable} it is
prosolvable relative to the diagonal representation $\rho$.  Thus
by the universal mapping property for
$\mathcal{S}_\rho(\pi_1(M,x))$ (see Section~\ref{relsolv}) there
exists an injection
\[
H^0 \left( EB(M,x)^L \right) \hookrightarrow \mathcal{O} \left(
\mathcal{S}_\rho \left( \pi_1(M, x) \right) \right).
\]
To show that this map is surjective it will suffice to show that
exponential iterated integrals compute the coordinate ring of any
solvable representation of $\p1Mx$ relative to $\rho$.
Our method is similar to that in Hain's proof of the 
$\pi_1$ de Rham theorem for ordinary iterated integrals in \cite{hain}.
Henceforth let $\pi = \pi_1(M,x)$. Suppose $\psi 
\colon \pi \to S$ is a solvable representation relative to $\rho$.  By
Corollary~\ref{relsolisupper}, we may assume that $S$ is a group
of upper triangular $n \times n$ matrices whose diagonal entries
are from $\mathcal{O}(T)$, meaning that they may be written as
$\int e^{\delta_k}$ with $\delta_k \in L$, $k=1, \ldots, n$.

Let $\CC^n = V^n \supset V^{n-1} \supset \ldots \supset V^0 = \{
0 \}$ denote the standard filtration,
\[
V^k = \left\{ (z_1, \ldots, z_k, 0, \ldots, 0) \mid z_1, \ldots, z_k \in
\CC \right\},
\]
which is stabilized by $S$.

Let $\widetilde{M}$ denote the universal cover of $M$.  From the
trivial flat bundle $\CC^n \times \widetilde{M}$ we can obtain a
bundle over $M$ with monodromy $\psi$: let
\[
E = \pi \backslash (\CC^n \times \widetilde{M} )
\]
where $\pi$ acts via $g \cdot (v, m) \mapsto \left( \psi(g)(v), mg^{-1}
\right)$.  The filtration $V^n \supseteq V^{n-1} \supseteq \ldots
\supseteq V^0$ induces a filtration of bundles
\[
E = E^n \supseteq E^{n-1} \supseteq \ldots \supseteq E^0 = 0.
\]
where the monodromy of the line bundle $E^k / E^{k-1}$ is given by
$\int e^{\delta_k}$.  We obtain for each $k$ a trivialization
\[
\CC \times M \to E^k / E^{k-1}
\]
via transport from $x$ with respect to the connection $d -
\delta_k$.  Composing these maps with splittings $E^k / E^{k-1}
\to E^k$, yields maps $\CC \times M \to E^k$ for $k = 1, \ldots,
n$. Adding these maps together we obtain an isomorphism
\[
\CC^n \times M \to E
\]
The induced connection form on $\CC^n \times M$ is an upper
triangular matrix with diagonal entries $\delta_1, \ldots,
\delta_n$.  By Proposition~\ref{upper}, the monodromy
representation $\rho \colon \pi \to S$ is equal to a matrix of
exponential iterated integrals with exponents from $\{ \delta_1,
\ldots, \delta_n\} \subseteq L$.  These matrix entries generate
the ring $\mc{O}(S)$ and this completes the proof.
\end{pf}

It is natural now to ask which diagonalizable representations of
$\pi_1(M,x)$ are defined by a module of closed $1$-forms.

\begin{proposition}
If $\rho \colon \p1Mx \rightarrow T \subseteq \left( \CC^*
\right)^n$ is a diagonalizable representation, then there exists a
defining $\ZZ$-module $L \subseteq E^1 (M; \CC)$ for $\rho$ if and
only if the induced map $H_1 \left( M \right) \rightarrow T$ is
trivial on $\Tor \left(H_1 \left( M \right) \right)$.
\end{proposition}

\begin{pf}
If $\delta$ is a closed $1$-form then the additive homomorphism
\[
\int \delta \colon H^1(M) \to \CC
\]
kills $\Tor H^1(M)$, and the same is true for $\int e^\delta
\colon H^1(M) \to \CC^*$.  So any diagonalizable representation
defined by $1$-forms kills $\Tor H^1(M)$.

For the converse, suppose $\rho \colon \pi \to \CC^*$ is a
homomorphism that kills $\Tor H^1(M)$.  Then $\rho$ induces a map
$\rho' \colon H_1(M; \CC) \to \CC^*$. Choose a basis $\{ z_i
\}_i$ for $H_1(M;\CC)$, and choose elements $f_i \in \CC$ such
that $e^{f_i} = \rho'(z_i)$.  By the ordinary de Rham theorem we
may find a closed $1$-form $\delta$ such that $\int \delta$ takes
$z_i$ to $f_i$.  Thus $\rho \cong \int e^\delta$.  This method
extends easily to define arbitrary diagonalizable representations
in terms of $1$-forms.
\end{pf}

Combining this result with the de Rham theorem gives a description
for
\begin{eqnarray}
\label{ring} H^0 \left( EB(M,x)^{B^1 \left(M; \CC \right)}
\right),
\end{eqnarray}
the ring of all closed exponential iterated integrals with closed
exponents.  If $\pi \to \textbf{G}$ is an algebraic representation
with Zariski dense image, then $\mathcal{O}(\textbf{G})$ is
computed by exponential iterated integrals with closed exponents
if and only if the reductive quotient $\textbf{G}/\textbf{G}_u$ is a 
diagonalizable
representation that kills $\Tor H^1(M)$.  The ring (\ref{ring})
is the direct limit of $\mathcal{O}(\textbf{G})$ over all such
representations.

\section{The Unipotent Radical of the Relative Solvable Completion}

\label{radical}

Recall that $H \to \mc{U}(H)$ denotes the unipotent completion of
the group $H$, or equivalently, the solvable completion of $H$
relative to the trivial representation $H \to \{ 1 \}$. The
functor $\mc{S}_\rho$ is right exact in the sense that
\[
\mc{U}(\ker \rho) \to \mc{S}_\rho(G) \to \mc{S}_\rho(\im \rho)
\to 1
\]
is exact for any group $G$ and diagonalizable representation
$\rho \colon G \to T$.  This is easily seen from the universal
mapping property.

The main theorem of this section asserts conditions under which
$\mc{U}(\ker \rho) \to \mc{S}_\rho(G)$ is an injection. The proof
will require an understanding of how $G$ acts by conjugation on
the solvable completions of its subgroups.  We begin with a
discussion of the automorphism groups of prounipotent and
relative prosolvable groups.

We will make free use of the equivalence of categories,
\[
\{\textnormal{prounipotent algebraic groups} / \CC \} \longleftrightarrow 
\{\textnormal{pronilpotent Lie algebras} / \CC \}
\]
provided by the {\bf exp} and {\bf log} maps.

Suppose $\mc{U}$ is a prounipotent group, and $\mc{U} = \mc{U}_0 \supseteq
\mc{U}_1 \supseteq \mc{U}_2 \supseteq \ldots$ is a filtration by closed
normal subgroups such that $\bigcap_{i = 0}^\infty \mc{U}_i = \{ 1\}$. Let
$\Aut_{\{ \mc{U}_i \}} \mc{U}$ denote the group of automorphisms of $\mc{U}$
that stabilize $\{ \mc{U}_i \}$.\footnote{To be precise,
$\Aut_{\{ \mc{U}_i \}} \mc{U}$ is the functor that takes a $\CC$-algebra 
$A$ to the group of automorphisms of $\mc{U} \times_\CC \Spec A$ that 
preserve $\{ \mc{U}_i \times_\CC \Spec A\}$.} Consider the graded vector 
space
obtained from the Lie algebra $\mc{u}$ of $\mc{U}$ with the induced
filtration $\{ \mc{u}_i \}$: \[ \gr_{\{\mc{u}_i\}} \mc{u} =
\mc{u}/\mc{u}_1 \oplus \mc{u}_1/\mc{u}_2 \oplus \mc{u}_2 / \mc{u}_3 \oplus
\ldots \]
The kernel of the morphism 
\begin{eqnarray}
\phi \colon
\Aut_{\{\mc{u}_i\}} \mc{u} \to \Aut \gr_{\{\mc{u}_i\}} \mc{u} 
\end{eqnarray} is 
prounipotent, and the same is true of the equivalent 
morphism
\begin{eqnarray}
\label{kerPhi}
\Phi \colon \Aut_{\{\mc{U}_i\}} \mc{U}  \to \Aut
\gr_{\{\mc{u}_i\}} \mc{u}.
\end{eqnarray}

If the chosen filtration is the central series $\{ \mc{U}^{(1)}, 
\mc{U}^{(2)}, \ldots
\}$ of $\mc{U}$, $\gr_{\{\mc{u}^{(i)}\}} \mc{u}$ has the
additional structure of a graded Lie algebra generated by $\mc{u}
/ \mc{u}^{(2)}$.  So an automorphism of $\gr_{\{\mc{u}^{(i)}\}}
\mc{u}$ is determined by its action on $\mc{u} / \mc{u}^{(2)}$:
\begin{eqnarray}
\label{Phi} \Phi \colon \Aut_{\{\mc{U}^{(i)}\}} \mc{U} \to \Aut
\gr_{\{\mc{u}^{(i)}\}} \mc{u} \hookrightarrow \GL ( \mc{u} /
\mc{u}^{(2)}).
\end{eqnarray}

Now suppose that 
\begin{eqnarray*}
\xymatrix{ 1 \ar[r] & \mc{U} \ar[r] & \mc{S} \ar[r]^\rho & T
\ar[r] & 1,}
\end{eqnarray*}
is an exact sequence where $T$ is diagonalizable.  Let
$\Aut_{\rho, \{\mc{U}_i\}} \mc{S}$ denote the group of
automorphisms of $\mc{S}$ that are $\rho$-invariant and stabilize
$\{\mc{U}_i\}$.  We show that the kernel of the morphism
\begin{eqnarray}
\label{solvmorph} \Aut_{\rho, \{\mc{U}_i\}} \mc{S} \to \Aut
\gr_{\{\mc{u}_i\}} \mc{u},
\end{eqnarray}
like $\ker \Phi$ above
(\ref{kerPhi}), is prounipotent. Fix (by
Proposition~\ref{itsplits}) a diagonalizable subgroup $T_0
\subseteq \mc{S}$ that maps isomorphically onto $T$.  Let $\ker
\Phi$ act on $\mc{S} = T_0 \ltimes \mc{U}$ leaving $T_0$ fixed,
and let $\mc{U}$ act on $\mc{S}$ by conjugation.  These actions
induce a morphism,
\begin{eqnarray*}
\ker \Phi \ltimes \mc{U} & \to & \Aut_{\rho, \{ \mc{U}_i \}}
\mc{S}, \\
(\psi, u) & \mapsto & \psi(u (\cdot) u^{-1}).
\end{eqnarray*}
We claim that the image of this morphism is exactly the kernel of
(\ref{solvmorph}).  Suppose that $\psi \in \Aut_{\rho,
\{\mc{U}_i\}} \mc{S}$ is an automorphism that fixes $\gr_{\{
\mc{u}_i\} } \mc{u}$. Since $\psi(T_0) \subseteq \mc{S}$ is
another closed subgroup that maps isomorphically onto $T$, we may
find $u \in \mc{U}$ such that $u \psi(T_0) u^{-1} = T_0$. Since $u
\psi ( \cdot ) u^{-1}$ is $\rho$-invariant, it therefore fixes
$T_0$; and since it also fixes $\gr_{\{ \mc{u}_i \}} \mc{u}$, we may
find it in the image of $\ker \Phi$.  This proves the claim.  And
indeed $\ker \Phi \ltimes \mc{U}$ is prounipotent since both $\ker \Phi$ 
and $\mc{U}$ are.

We summarize this discussion.

\begin{lemma}
\label{AutS} Suppose that
\begin{eqnarray*}
\xymatrix{ 1 \ar[r] & \mc{U} \ar[r] & \mc{S} \ar[r]^\rho & T
\ar[r] & 1}
\end{eqnarray*}
is an exact sequence of groups in which $T$ is
diagonalizable and $\mc{U}$ is prounipotent. Suppose that $\mc{U}
= \mc{U}_0 \supseteq \mc{U}_1 \supseteq \mc{U}_2 \supseteq
\ldots$ is a filtration by closed normal subgroups such that
$\bigcap_{i = 0}^\infty \mc{U}_i = \{ 1\}$. Then the kernel of
the morphism
\begin{eqnarray*}
\Aut_{\rho, \{\mc{U}_i\}} \mc{S} \to \Aut \gr_{\{\mc{u}_i\}}
\mc{u}
\end{eqnarray*}
is prounipotent.
\end{lemma}

We are now ready to state the main theorem.  Suppose that $1 \to K \to G
\to A \to 1$ is an exact sequence of groups with $A$ abelian and finitely
generated.  Note that on the vector space $H_1(K; \CC)$ (or equivalently,
$\mc{U}(K) / [\mc{U}(K), \mc{U}(K)]$) there is an induced conjugation
action $G \to \Aut H_1(K; \CC)$, which is abelian. We are interested in
the case when this action is algebraic, and thus may be written as \[ G
\to T \times \left( \GG_a \right)^m \to \Aut H_1(K; \CC), \] where $T$ is
a diagonalizable group.

\begin{theorem}
\label{unipotentinjects}
Suppose $1 \to K \to G \to A \to 1$ is an exact sequence of groups with 
$A$ abelian and finitely generated, and suppose that the 
conjugation action $G \to \Aut H_1(K;\CC)$ factors as
\[
G \to T \times \left( \GG_a \right)^m \to \Aut H_1(K; \CC)
\]
where $\rho \colon G \to T$ is a diagonalizable representation with 
Zariski dense image.  Then 
\begin{eqnarray}
\label{unipexact}
1 \to \mc{U}(K) \to \mc{S}_\rho(G) \to \mc{S}_\rho(\im \rho) \to 1
\end{eqnarray}
is an exact sequence.
\end{theorem}

\begin{pf}
Let
\[
A = \ZZ / (n_1) \oplus \ZZ / (n_2) \oplus \ldots \oplus \ZZ /
(n_r),
\]
with $n_i \in \ZZ$.  We induct on $r$. Let $G' \subseteq G$ be
the inverse image of $\ZZ / (n_1) \oplus \ldots \oplus \ZZ /
(n_{r-1})$, and assume that $\mc{U}(K) \hookrightarrow
\mc{S}_{\rho} (G')$.

Let $\overline{G}$ denote the pushout of the diagram
\[
\xymatrix{G' \ar[r] \ar[d] & G \\
\mc{S}_\rho(G')}.
\]
Then $\overline{G}$ fits into an exact sequence
\begin{eqnarray}
\label{oftheseq} 1 \to \mc{S}_{\rho}(G') \to \overline{G} \to \ZZ
/ (n_r) \to 0.
\end{eqnarray}
For convenience, we consider $\mc{U}(K)$ and $\mc{S}_{\rho}(G')$ as
subgroups of $\overline{G}$.

We desire a splitting of the
sequence (\ref{oftheseq}). If $n_r =
0$ this is straightforward. Otherwise, choose an element of $G$
that maps to $(0, 0, \ldots, 0, 1) \in A$, and let $g$ denote its
image in $\overline{G}$. Then $g^{n_r}$ is contained in the
prounipotent group $\mc{U}(K)$, so it has a unique $n_r$th root
$u$ in $\mc{U}(K)$, which commutes with $g$. Replacing $g$ with
$u^{-1}g$ we have $g^{n_r} = 1$, hence
\[
\overline{G} = \left< g \right> \ltimes \mc{S}_\rho(G').
\]
Lemma~\ref{AutS} will show that the
conjugation action of $g$ on $\mc{S}_{\rho}(G')$ is algebraic.
Filter the prounipotent kernel of $\mc{S}_{\rho}(G') \to T$ via
the central series of $\mc{U}(K)$:
\begin{eqnarray*}
\mc{U}_0 & = & \ker [\mc{S}_{\rho}(G') \to T] \\
\mc{U}_1 & = & \mc{U}(K) \\
\mc{U}_2 & = & [\mc{U}_1, \mc{U}_1] \\
\mc{U}_3 & = & [\mc{U}_1, \mc{U}_2] \\
\mc{U}_4 & = & [\mc{U}_1, \mc{U}_3] \\
& \ldots
\end{eqnarray*}
This filtration is evidently stabilized by $g$.  The action of $g$ on 
$\mc{u}_0/\mc{u}_1$ is trivial, and the action of $g$ on the graded Lie 
algebra $\gr_{\{ \mc{u}_1, \mc{u}_2, \ldots  \}} \mc{u}$ is isomorphic to 
the action of $g$ on $\mc{u}_1 / \mc{u}_2 \cong H_1(\mc{U}; \CC)$  
(recall the discussion prior to (\ref{Phi})), hence the 
commutative diagram
\[
\xymatrix{ & & \Aut_{\rho, \left\{ \mc{U}_i \right\} } \mc{S}_\rho(G') 
\ar[d] \\
\left< g \right> \ar[r] \ar[rru] & T \times \left( \GG_a \right)^m \ar[r] 
&
\gr_{\{ \mc{u}_i \}} \mc{u}.
}
\]
The horizontal map $\left< g \right> \to \gr_{\{ \mc{u}_i \}} \mc{u}$
factors through $\mc{S}_\rho( \left< g \right> )$, and the downward map 
has
prounipotent kernel (Lemma~\ref{AutS}), therefore the diagonal map also 
factors
through $\mc{S}_\rho ( \left< g \right> )$ by the universal mapping 
property.  Thus we can ``thicken'' $\overline{G}$ further:
\[
\overline{G} = \left< g \right> \ltimes \mc{S}_{\rho}(G')
\hookrightarrow \mc{S}_\rho(\left< g\right> ) \ltimes
\mc{S}_{\rho}(G').
\]
This completes the proof, as $\mc{S}_\rho(\left< 
g\right> ) \ltimes \mc{S}_{\rho}(G')$ contains $\mc{U}(K)$ and
is isomorphic to
$\mc{S}_\rho(G)$, as can be seen from the universal mapping
property.
\end{pf}

\section{Exponential Iterated Integrals on Complex Curves and Fibered 
Knots}

\label{examples}

In this section we consider two examples that illustrate the $\pi_1$ de
Rham theorems and Theorem~\ref{unipotentinjects}.

We are interested in
manifolds for which the ring of closed exponential iterated integrals can
be calculated explicitly.  We have earlier referred to the problem that
the same map $P_{x,x} M \to \CC$ may be computed by two different integral
expressions.  This problem can sometimes be solved by putting a complex
structure on $M$ and restricting attention to integrals composed from
holomorphic $1$-forms.

Consider 
first the case when $M$ is a smooth affine complex curve.  Since $\p1Mx$ 
is free, the unipotent completion $\p1Mx \to \mc{U}(\p1Mx )$ is a faithful 
representation,\footnote{For a proof see Appendix A3 in \cite{quillen}, 
Proposition 3.6(a) in particular.} hence ordinary iterated integrals are 
adequate to distinguish any two elements of $\p1Mx$.  The next proposition
provides a description for $H^0(B(M,x))$.

\begin{proposition}
\label{thebasis}
Suppose $(M,x)$ is a smooth affine complex curve.  Closed iterated 
integrals separate the elements of $\p1Mx$.  
Let $\{ \omega_1, \ldots, 
\omega_n \}$ be a basis for the space of closed holomorphic $1$-forms on 
$M$.  Then
\begin{eqnarray}
\mc{B} = \{ 1 \} \cup \left\{ \int \omega_{j_1} \ldots \omega_{j_k} \mid k 
> 0, (j_i) \in \left\{ 1, \ldots, n \right\}^k \right\}
\end{eqnarray}
is a basis for $H^0(B(M,x))$.
\end{proposition}

\begin{pf}
The key observations are that each class in $H^1(M; \CC)$ has a unique 
holomorphic representative, and that the wedge product of any two closed 
holomorphic $1$-forms is zero.  The fact that $\mc{B}$ is a basis 
for $H^0(B(M,x))$ actually 
follows from a 
general Theorem of Chen's (\cite{chen2}, Theorem 4.1.1) which deals with 
the complex of higher iterated integrals.  We provide here an elementary 
proof for the sake of clarity.

\begin{lemma}
Any connection on $M$ of the form $d - \theta$, where $\theta$ is a matrix 
of closed holomorphic $1$-forms, is flat.
Any flat connection $d
- \omega$, with $\omega$ a nilpotent upper triangular matrix of $1$-forms,
is conjugate via a matrix $G \colon M \to \UnC$ to such a connection $d - 
\theta$.
\end{lemma}

\begin{pf}
The first assertion is immediate since the curvature of $d - \theta$ is $d 
\theta + \theta \wedge \theta = 0$.  For the second, take any $k \in 
\{1, \ldots, 
n-1\}$ and suppose
\begin{eqnarray}
d - \left[ \begin{array}{cccccc}
0 & \omega_{12} & \omega_{13} & \omega_{14} & \cdots & \omega_{1n} \\
0 & 0 & \omega_{23} & \omega_{24} & \cdots & \omega_{2n} \\
0 & 0 & 0 & \omega_{34} & \cdots & \omega_{3n} \\
\vdots & \vdots & \vdots & \vdots & \ddots & \vdots \\
0 & 0 & 0 & 0 & \cdots & 0 \\
\end{array} \right].
\end{eqnarray}
is a flat connection such that each entry $\omega_{ij}$ with $j - i < k$ 
is closed and holomorphic.  The equation $d \omega + \omega 
\wedge \omega = 0$ implies that the $1$-forms $\omega_{ij}$ for which $j - 
i = k$ must be closed.  Choose for each such entry an exact 
$1$-form $df_{ij}$ such that $\omega_{ij} - df_{ij}$ is holomorphic. Let 
$F$ be the $n \times n$ matrix with $F_{ij}$ equal to $f_{ij}$ if $j - i = k$ 
and zero otherwise.  The matrix $\omega'$ satisfying 
\begin{eqnarray}
\label{trans}
d - \omega' = (I + F)^{-1}(d - \omega)(I + F)
\end{eqnarray}
has $\omega'_{ij} = \omega_{ij} - dF_{ij}$ for $j - i \leq k$.  
Continuing 
by induction on $k$ we obtain the desired matrix $\theta$.
\end{pf}

Now since any closed iterated integral $I$ on $M$ arises from the 
monodromy 
of such a connection $d - \omega$ (recall the proof of 
Theorem~\ref{mainthm}), it may be expressed in terms of closed 
holomorphic $1$-forms.  And linearity
\[
\int \gamma_1 \ldots (z\gamma_i + z' \gamma'_i ) \ldots \gamma_n = z \int 
\gamma_1 \ldots \gamma_i \ldots \gamma_n + z' \int \gamma_1 \ldots 
\gamma'_i \ldots \gamma_n,
\]
allows us then to express $I$ as a sum of elements from $\mc{B}$.

It remains only to show that $\mc{B}$ is a linearly indepedent set;
this follows from a straightforward manipulation.  Let $\lambda_1, 
\ldots, \lambda_{n}$ be free generators for $\pi_1(M,x)$; it suffices to 
prove the proposition for any chosen basis 
$\{ \omega_j \}_j$, so let us assume $\{ \omega_j \}_j$ is such that $\{ 
[\omega_j] \}_j \subseteq H^1(M; \CC)$ is the dual basis for $\{ 
[\lambda_j ] \}_j \subseteq H_1(M; \CC)$.  
Suppose
\begin{eqnarray}
\label{sum}
\sum_{I \in \mc{B}' \subseteq \mc{B}} z_I I = 0
\end{eqnarray}
is a nontrivial finite sum with each $z_I$ nonzero.  Evidently such 
an expression must include a term of length at least $2$.  Suppose
$z_{I_0} I_0 = z_{I_0} \int \omega_{i_1} \ldots \omega_{i_k}$ 
is a term of maximal length.  The transformation
\[
I \mapsto  I - \left< \Delta I, ( \cdot, \lambda_{i_k} ) \right> 
\]
(applying the comultiplication formula (\ref{comultdef}) formally)
turns (\ref{sum}) into a nontrivial sum with terms of 
length $\leq k - 1$.  Continuing in this manner yields a contradiction.
This completes the proof.
\end{pf}

Now we consider the example of a fibered knot complement.  Suppose $K 
\subseteq S^3$ is a tame knot and that $S^3 \smallsetminus K$ has an 
infinite-cyclic covering map
\[
\phi \colon (\RR \times F, (0, \overline{x})) \to (S^3 \smallsetminus K, 
x),
\]
whose deck transformations $\Psi_n, n \in \ZZ$ are given by 
$\Psi_1(t,f) = (t+1, \psi(f))$ where $\psi \colon F \to F$ is a 
homeomorphism.  $F$ is a noncompact $2$-manifold which we may take to 
be an affine complex curve, making the previous discussion useful.
Note that pulling back an 
exponential iterated integral on $S^3 \smallsetminus K$ with closed exponents,
\[
\phi^* \int e^{\delta_0} \omega_{01} \ldots \omega_{(n-1)n} e^{\delta_n} = 
\int e^{\phi^* \delta_0} \phi^* \omega_{01} \ldots \phi^* \omega_{(n-1)n} 
e^{\phi^* \delta_n},
\]
produces an exponential iterated integral with exact exponents, which may 
be rewritten via Proposition~\ref{expexact} as an ordinary iterated 
integral.  Let $L$ be the $\ZZ$-module of closed $1$-forms that defines 
the diagonal part of the conjugation representation
\[
\pi_1(S^3 \smallsetminus K ,x) \to \Aut H_1(\RR \times F; \CC ).
\]
The content of Theorem~\ref{unipotentinjects} is that the map
\begin{eqnarray}
\label{phistar}
\phi^* \colon H^0(EB(S^3 \smallsetminus K,x)^L) \to H^0(B(\RR \times F, 
(0, \overline{x}))) = H^0(B(F, \overline{x}))
\end{eqnarray}
is surjective.  Combining this fact with Proposition~\ref{thebasis} shows 
that $H^0(EB(S^3 \smallsetminus K, x)^L)$ separates the elements of 
$\pi_1(S^3 \smallsetminus K, x)$.  It also allows one to describe 
$H^0(EB(S^3 \smallsetminus K, x)^L)$ by looking at the pre-image of the 
basis $\mc{B}$.  We demonstrate this approach with the 
particular example of the trefoil knot.

It is convenient (for coordinates) to consider a complex manifold that is
of the same homotopy type as the complement of the trefoil.  Assume \[ S^3
= \{ (x,y) \mid \left| x \right|^2 + \left| y \right|^2 = 1 \}
\subseteq \CC^2. \] The knot,
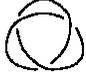
\[
\xygraph{!{0;/r.75pc/:} !P3"a"{~>{}} !P12"b"{~:{(1.414,0):}~>{}}
!{\vover~{"b2"}{"b1"}{"a1"}{"a3"}} !{\save 0;"b2"-"b5":"b5",
  \xcaph @(+.1)\restore}
!{\vover~{"b6"}{"b5"}{"a2"}{"a1"}} !{\save 0;"b6"-"b9":"b9",
  \xcaph @(+.2)\restore}
!{\vover~{"b10"}{"b9"}{"a3"}{"a2"}} !{\save 0;"b10"-"b1":"b1",
   \xcaph @(+.3)\restore} }
\]
can be specified as the 
intersection of $S^3$ with the singular curve
\[
C = \{ (x, y) \mid x^3 + y^2 = 0 \},
\]
and $\CC^2 \smallsetminus C$ deformations retracts onto $S^3 
\smallsetminus C \cap S^3$.  The infinite cyclic cover of $\CC^2 
\smallsetminus C$ may written as $\CC \times F$, where $F = \{ (x, y) 
\mid x^3 + y^2 = 1 \} \subseteq \CC^2$:
\begin{eqnarray*}
\phi \colon \CC \times F & \to & \CC^2 \smallsetminus C \\
(t, x, y) & \mapsto & (e^{\frac{2 \pi i}{3}t} x, e^{\pi i t}y).
\end{eqnarray*}
The homeomorphism $\psi \colon F \to F$ that generates the group of deck 
transformations is given by
\[
\psi(x, y) = (\zeta_6^{-2} x, \zeta_6^{-3} y).
\]
The closed holomorphic $1$-forms
\[
\omega_{-1} = \frac{dx}{y}, \hskip0.25in \omega_1 = \frac{xdx}{y}
\]
on $F$ diagonalize the action of $\psi^*$ on $H^1(F; \CC)$, for
\[
\psi^* \omega_{-1} = \zeta_6 \omega_{-1}, \hskip0.25in \psi^* \omega_1 = 
\zeta_6^{-1} 
\omega_1.
\]
Extend the $1$-forms $\omega_{-1}$ and 
$\omega_1$ to $\CC \times F$ by pulling back along 
the projection map to $F$.  The $1$-forms
\[
e^{-\frac{\pi i}{3}t} \omega_{-1}, e^{\frac{\pi i}{3}t} \omega_1, 
\frac{\pi i}{3} dt \in 
E^1(\CC \times F ; \CC)
\]
are each invariant under the action of the group of deck transformations.
Let $\omega_1^\circ$, $\omega_{-1}^\circ$, and $\delta$ denote 
the $1$-forms that they induce on $\CC^2 \smallsetminus C$, respectively.

Consider the integral
\begin{eqnarray}
\label{oftheform}
I =
\int \omega_{\epsilon_1}^\circ e^{\epsilon_1 \delta} 
\omega_{\epsilon_2}^\circ
e^{(\epsilon_1 + \epsilon_2)\delta} \ldots \omega_{\epsilon_n}^\circ 
e^{(\epsilon_1 + \ldots + \epsilon_n)\delta},
\end{eqnarray}
where $\epsilon_1, \ldots, \epsilon_n \in \{ \pm 1 \}$.  Applying 
Proposition~\ref{expexact} we find
\begin{eqnarray*}
\phi^* I 
& = & \int \omega_{\epsilon_1} \omega_{\epsilon_2} \ldots 
\omega_{\epsilon_n}. \\
\end{eqnarray*}
Integrals of this form provide a basis for $H^0(B(\CC \times F, (0, 
\overline{x})))$, as we know.
The exact sequence (\ref{unipexact}) from 
Theorem~\ref{unipotentinjects} implies that $H^0(EB(\CC^2 
\smallsetminus C, x)^{\left< \delta \right>})$ is isomorphic as a 
$\CC$-algebra to $H^0(B(\CC \times F, (0, 
\overline{x}))) \otimes \mc{O}(\mc{S}_\rho (\im \rho))$, hence the 
following.

\begin{proposition}
The set of integrals 
\begin{eqnarray}
\label{theel}
\int \delta^m \int e^{k \delta} 
\int \omega_{\epsilon_1}^\circ e^{\epsilon_1 \delta} 
\omega_{\epsilon_2}^\circ
e^{(\epsilon_1 + \epsilon_2)\delta} \ldots \omega_n^\circ e^{(\epsilon_1 + 
\ldots + \epsilon_n) \delta}
\end{eqnarray}
with $k \in \ZZ$, $m, n \geq 0$, and $\epsilon_1, \ldots, \epsilon_n \in 
\{ \pm 1 \}$, is a basis for $H^0(EB(\CC^2 \smallsetminus C, x)^{\left< 
\delta \right>})$.
\end{proposition}

The comultiplication formula (\ref{comultdef}) is easily applied to these 
integrals, giving us an effective description for the Hopf
algebra $H^0(EB(\CC^2 \smallsetminus C, x)^{\left< \delta \right>})$.

\section{Acknowledgements}

This research was started in May 2000 at PRUV, a summer program at Duke
University supported by a VIGRE grant.  Many thanks to my undergraduate
advisor Richard Hain for his generous assistance and guidance.  The idea
and outline for the paper were his.  Thanks also to Arthur Ogus and
Chee-Whye Chin for helpful discussions pertaining to 
Section~\ref{radical}.



\end{document}